\documentclass[12pt]{article}

\usepackage{amsmath,amsthm}
\usepackage{amssymb}
\usepackage{color}
\usepackage{xspace}
\usepackage[colorlinks=true,
linkcolor=green,
filecolor=brown,
citecolor=green]{hyperref}

\def\s{r}
\def\M{M}

\catcode`\@=11

\thicklines
\newskip\Einheit \Einheit=.6cm
\newcount\xcoord \newcount\ycoord
\newdimen\xdim \newdimen\ydim \newdimen\PfadD@cke \newdimen\Pfadd@cke
\def\PfadDicke#1{\PfadD@cke#1 \divide\PfadD@cke by2 
\Pfadd@cke\PfadD@cke \multiply\PfadD@cke by2}
\long\def\LOOP#1\REPEAT{\def\BODY{#1}\ITERATE}
\def\ITERATE{\BODY \let\next\ITERATE \else\let\next\relax\fi \next}
\let\REPEAT=\fi
\def\Punkt{\hbox{\raise-2pt\hbox to0pt{\hss\scriptsize$\bullet$\hss}}}

\def\DuennPunkt(#1,#2){\unskip
  \raise#2 \Einheit\hbox to0pt{\hskip#1 \Einheit
          \raise-1.5pt\hbox to0pt{\hss\tiny$\bullet$\hss}\hss}}
		  
\def\NormalPunkt(#1,#2){\unskip
  \raise#2 \Einheit\hbox to0pt{\hskip#1 \Einheit
          \raise-3pt\hbox to0pt{\hss\large$\bullet$\hss}\hss}}
\def\DickPunkt(#1,#2){\unskip
  \raise#2 \Einheit\hbox to0pt{\hskip#1 \Einheit
          \raise-4pt\hbox to0pt{\hss\Large$\bullet$\hss}\hss}}
\def\Kreis(#1,#2){\unskip
  \raise#2 \Einheit\hbox to0pt{\hskip#1 \Einheit
          \raise-4pt\hbox to0pt{\hss\Large$\circ$\hss}\hss}}
\def\Diagonale(#1,#2)#3{\unskip\leavevmode
  \xcoord#1\relax \ycoord#2\relax
      \raise\ycoord \Einheit\hbox to0pt{\hskip\xcoord \Einheit
         \unitlength\Einheit
         \line(1,1){#3}\hss}}
\def\AntiDiagonale(#1,#2)#3{\unskip\leavevmode
  \xcoord#1\relax \ycoord#2\relax \advance\xcoord by -0.05\relax
      \raise\ycoord \Einheit\hbox to0pt{\hskip\xcoord \Einheit
         \unitlength\Einheit
         \line(1,-1){#3}\hss}}
\def\Pfad(#1,#2),#3\endPfad{\unskip\leavevmode
  \xcoord#1 \ycoord#2 \thicklines\ZeichnePfad#3\endPfad\thinlines}
\def\ZeichnePfad#1{\ifx#1\endPfad\let\next\relax
  \else\let\next\ZeichnePfad
    \ifnum#1=1
      \raise\ycoord \Einheit\hbox to0pt{\hskip\xcoord \Einheit
         \vrule height\Pfadd@cke width1 \Einheit depth\Pfadd@cke\hss}%
      \advance\xcoord by 1
     \else\ifnum#1=2
      \raise\ycoord \Einheit\hbox to0pt{\hskip\xcoord \Einheit
         \unitlength\Einheit
         \line(0,1){1}\hss}
      \advance\xcoord by 0
      \advance\ycoord by 1
 \else\ifnum#1=3
      \raise\ycoord \Einheit\hbox to0pt{\hskip\xcoord \Einheit
         \unitlength\Einheit
         \line(1,1){1}\hss}
      \advance\xcoord by 1
      \advance\ycoord by 1
    \else\ifnum#1=4
      \raise\ycoord \Einheit\hbox to0pt{\hskip\xcoord \Einheit
         \unitlength\Einheit
         \line(1,-1){1}\hss}
      \advance\xcoord by 1
      \advance\ycoord by -1
   \else\ifnum#1=5
      \raise\ycoord \Einheit\hbox to0pt{\hskip\xcoord \Einheit
         \unitlength\Einheit
         \line(2,1){2}\hss}
      \advance\xcoord by 2
      \advance\ycoord by 1
	  \else\ifnum#1=6
      \raise\ycoord \Einheit\hbox to0pt{\hskip\xcoord \Einheit
         \unitlength\Einheit
         \line(2,-1){2}\hss}
      \advance\xcoord by 2
      \advance\ycoord by -1
	  \else\ifnum#1=7
      \raise\ycoord \Einheit\hbox to0pt{\hskip\xcoord \Einheit
         \unitlength\Einheit
         \line(3,1){3}\hss}
      \advance\xcoord by 3
      \advance\ycoord by 1
	  \else\ifnum#1=8
      \raise\ycoord \Einheit\hbox to0pt{\hskip\xcoord \Einheit
         \unitlength\Einheit
         \line(3,-1){3}\hss}
      \advance\xcoord by 3
      \advance\ycoord by -1
    \fi\fi\fi\fi\fi\fi\fi\fi
  \fi\next}
\def\hSSchritt{\leavevmode\raise-.4pt\hbox 
to0pt{\hss.\hss}\hskip.2\Einheit
  \raise-.4pt\hbox to0pt{\hss.\hss}\hskip.2\Einheit
  \raise-.4pt\hbox to0pt{\hss.\hss}\hskip.2\Einheit
  \raise-.4pt\hbox to0pt{\hss.\hss}\hskip.2\Einheit
  \raise-.4pt\hbox to0pt{\hss.\hss}\hskip.2\Einheit}
\def\vSSchritt{\vbox{\baselineskip.2\Einheit\lineskiplimit0pt
\hbox{.}\hbox{.}\hbox{.}\hbox{.}\hbox{.}}}
\def\DSSchritt{\leavevmode\raise-.4pt\hbox to0pt{%
  \hbox to0pt{\hss.\hss}\hskip.2\Einheit
  \raise.2\Einheit\hbox to0pt{\hss.\hss}\hskip.2\Einheit
  \raise.4\Einheit\hbox to0pt{\hss.\hss}\hskip.2\Einheit
  \raise.6\Einheit\hbox to0pt{\hss.\hss}\hskip.2\Einheit
  \raise.8\Einheit\hbox to0pt{\hss.\hss}\hss}}
\def\dSSchritt{\leavevmode\raise-.4pt\hbox to0pt{%
  \hbox to0pt{\hss.\hss}\hskip.2\Einheit
  \raise-.2\Einheit\hbox to0pt{\hss.\hss}\hskip.2\Einheit
  \raise-.4\Einheit\hbox to0pt{\hss.\hss}\hskip.2\Einheit
  \raise-.6\Einheit\hbox to0pt{\hss.\hss}\hskip.2\Einheit
  \raise-.8\Einheit\hbox to0pt{\hss.\hss}\hss}}
\def\SPfad(#1,#2),#3\endSPfad{\unskip\leavevmode
  \xcoord#1 \ycoord#2 \ZeichneSPfad#3\endSPfad}
\def\ZeichneSPfad#1{\ifx#1\endSPfad\let\next\relax
  \else\let\next\ZeichneSPfad
    \ifnum#1=1
      \raise\ycoord \Einheit\hbox to0pt{\hskip\xcoord \Einheit
         \hSSchritt\hss}%
      \advance\xcoord by 1
    \else\ifnum#1=2
      \raise\ycoord \Einheit\hbox to0pt{\hskip\xcoord \Einheit
        \hbox{\hskip-2pt \vSSchritt}\hss}%
      \advance\ycoord by 1
    \else\ifnum#1=3
      \raise\ycoord \Einheit\hbox to0pt{\hskip\xcoord \Einheit
         \DSSchritt\hss}
      \advance\xcoord by 1
      \advance\ycoord by 1
    \else\ifnum#1=4
      \raise\ycoord \Einheit\hbox to0pt{\hskip\xcoord \Einheit
         \dSSchritt\hss}
      \advance\xcoord by 1
      \advance\ycoord by -1
    \fi\fi\fi\fi
  \fi\next}
\def\Koordinatenachsen(#1,#2){\unskip
 \hbox to0pt{\hskip-.5pt\vrule height#2 \Einheit width.5pt depth1 
\Einheit}%
 \hbox to0pt{\hskip-1 \Einheit \xcoord#1 \advance\xcoord by1
    \vrule height0.25pt width\xcoord \Einheit depth0.25pt\hss}}
\def\Koordinatenachsen(#1,#2)(#3,#4){\unskip
 \hbox to0pt{\hskip-.5pt \ycoord-#4 \advance\ycoord by1
    \vrule height#2 \Einheit width.5pt depth\ycoord \Einheit}%
 \hbox to0pt{\hskip-1 \Einheit \hskip#3\Einheit 
    \xcoord#1 \advance\xcoord by1 \advance\xcoord by-#3 
    \vrule height0.25pt width\xcoord \Einheit depth0.25pt\hss}}
\def\Gitter(#1,#2){\unskip \xcoord0 \ycoord0 \leavevmode
  \LOOP\ifnum\ycoord<#2
    \loop\ifnum\xcoord<#1
      \raise\ycoord \Einheit\hbox to0pt{\hskip\xcoord 
\Einheit\Punkt\hss}%
      \advance\xcoord by1
    \repeat
    \xcoord0
    \advance\ycoord by1
  \REPEAT}
\def\Gitter(#1,#2)(#3,#4){\unskip \xcoord#3 \ycoord#4 \leavevmode
  \LOOP\ifnum\ycoord<#2
    \loop\ifnum\xcoord<#1
      \raise\ycoord \Einheit\hbox to0pt{\hskip\xcoord 
\Einheit\Punkt\hss}%
      \advance\xcoord by1
    \repeat
    \xcoord#3
    \advance\ycoord by1
  \REPEAT}
\def\Label#1#2(#3,#4){\unskip \xdim#3 \Einheit \ydim#4 \Einheit
  \def\lo{\advance\xdim by-.5 \Einheit \advance\ydim by.5 \Einheit}%
  \def\llo{\advance\xdim by-.25cm \advance\ydim by.5 \Einheit}%
  \def\loo{\advance\xdim by-.5 \Einheit \advance\ydim by.25cm}%
  \def\o{\advance\ydim by.25cm}%
  \def\ro{\advance\xdim by.5 \Einheit \advance\ydim by.5 \Einheit}%
  \def\rro{\advance\xdim by.25cm \advance\ydim by.5 \Einheit}%
  \def\roo{\advance\xdim by.5 \Einheit \advance\ydim by.25cm}%
  \def\l{\advance\xdim by-.30cm}%
  \def\r{\advance\xdim by.30cm}%
  \def\lu{\advance\xdim by-.5 \Einheit \advance\ydim by-.6 \Einheit}%
  \def\llu{\advance\xdim by-.25cm \advance\ydim by-.6 \Einheit}%
  \def\luu{\advance\xdim by-.5 \Einheit \advance\ydim by-.30cm}%
  \def\u{\advance\ydim by-.30cm}%
  \def\ru{\advance\xdim by.5 \Einheit \advance\ydim by-.6 \Einheit}%
  \def\rru{\advance\xdim by.25cm \advance\ydim by-.6 \Einheit}%
  \def\ruu{\advance\xdim by.5 \Einheit \advance\ydim by-.30cm}%
  #1\raise\ydim\hbox to0pt{\hskip\xdim
     \vbox to0pt{\vss\hbox to0pt{\hss$#2$\hss}\vss}\hss}%
}
\catcode`\@=12

\begin{document}
\newtheorem{theorem}{Theorem}
\newtheorem{defn}[theorem]{Definition}
\newtheorem{lemma}[theorem]{Lemma}
\newtheorem{prop}[theorem]{Proposition}
\newtheorem{cor}[theorem]{Corollary}
\begin{center}
{\Large
An application of a bijection of Mansour, Deng, and Du                        \\ 
}

\vspace{10mm}
David Callan \\ 
\vspace*{-1mm}
Department of Statistics,  University of Wisconsin-Madison  \\
\vspace*{-1mm}
1300 University Avenue, Madison, WI \ 53706-1532  \\
{\bf callan@stat.wisc.edu}  \\
\vspace{5mm}

October 24, 2012

\end{center}

\begin{abstract}
The large Schr\"{o}der numbers are known to count several classes of permutations avoiding two 4-letter patterns. Here we show they count another family of permutations, 
those whose left to right minima decomposition, when reversed, is 321-avoiding. The main tool is the Mansour-Deng-Du bijection from 321-avoiding permutations to Dyck paths.
\end{abstract}

\vspace{5mm}

{\Large \textbf{1 \quad Introduction}  }
The large Schr\"{o}der numbers $(\s_{n})$,  sequence \htmladdnormallink{A006318}{http://oeis.org/A006318}
in the OEIS \cite{oeis}, are well known to count Schr\"{o}der $n$-paths---nonnegative paths of upsteps $U=(1,1)$, downsteps $D=(1,-1)$, and double flatsteps $F=(2,0)$ from the origin $(0,0)$ to $(2n,0)$. Among their other combinatorial interpretations are several involving pattern avoidance in permutations, including separable permutations (i.e., those that avoid the patterns 2413 and 3142) and permutations sortable by an output-restricted deque (equivalently, those that avoid the patterns 2431 and 4231) \cite[Ex. 6.39\,(l,m)]{ec2}. Here we show they count another family of permutations.

The left to right minima decomposition of a permutation (in one-line notation) is obtained by splitting it just before each left to right minimum. Thus $\tau=4\,6\,5\,2\,3\,8\,1\,7$ decomposes as 4\,6\,5,\ 2\,3\,8,\ 1\,7. For a permutation $\pi$, we will denote by $f(\pi)$ the permutation obtained by reversing this list of subpermutations and concatenating. Thus $f(\tau)= 1\,7\,2\,3\,8\,4\,6\,5$. Our result is that $\s_{n}$ is the number of permutations $\pi$ of $[n+1]$ for which $f(\pi)$ is 321-avoiding. Of course, $f(\tau)$ fails to be  321-avoiding due to the 865. This family is not closed under containment---consider 3254 in 13254---and so is not a pattern avoidance class.

\newpage


{\Large \textbf{2 \quad Outline of proof}  } For a permutation $\pi$ of $[n+1]$, the first entry of $f(\pi)$ is always 1; let $f'(\pi)$ denote the result of deleting this initial 1 and decrementing all other entries by 1. Clearly, $f'(\pi)$ is a permutation of $[n]$ and is 321-avoiding if and only if $f(\pi)$ is.

Now let $\M$ denote the Mansour-Deng-Du bijection \cite{mansour2006} (see also \cite{revisit07}\,) from 321-avoiding permutations of $[n]$ to Dyck $n$-paths (a Dyck $n$-path is a Schr\"{o}der $n$-path with no flatsteps). This bijection is reviewed in somewhat simplified form in the next section. 
The key to the proof is the fact, noted in \cite{mansour2006}, that $\M$ takes the \emph{right to left} minima in a 321-avoiding permutation to the \emph{peaks} in the corresponding Dyck path, even preserving locations. So, given a permutation $\pi$ of $[n+1]$ for which $f(\pi)$ is 321-avoiding, apply $\M$ to $f'(\pi)$ to obtain a Dyck $n$-path $P$. Since, from the definitions of $f$ and $f'$, each left to right minimum of $\pi$ other than 1 is, after decrementing, a right to left minimum of $f'(\pi)$, it corresponds to a peak in $P$. Change this peak, $UD$, to a double flatstep $F$. As verified in the next section, this mapping is a bijection from the permutations of $[n+1]$ being counted to Schr\"{o}der $n$-paths, establishing the desired count.

\vspace*{10mm}

{\Large \textbf{3 \quad The Mansour-Deng-Du bijection}  }

First, we recall that the \emph{ascent-descent code} \cite[p.\,3]{revisit07} of a Dyck path (see also \cite[Ex. 6.19, item ($\textrm{i}^{6})$]{ec2}\,)
is obtained by recording all but the last of the partial sums of the ascent lengths (resp. descent lengths). For the path shown below, the ascent lengths are (1,3,3,1,2,1), the  descent lengths are (1,1,4,1,2,2) and so the ascent-descent code is 
$\left(\begin{smallmatrix}1&4&7&8&10 \\
1&2&6&7&9 \end{smallmatrix}\right)$. Each ascent ends with a peak upstep.

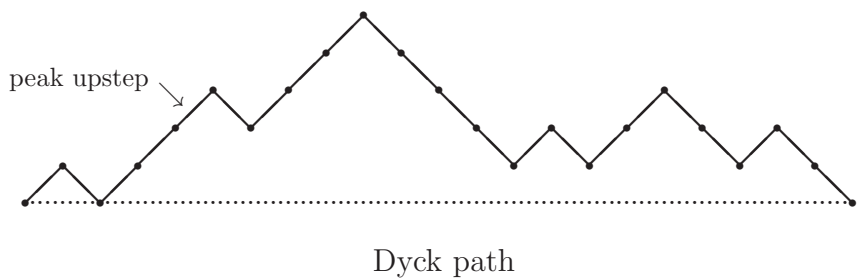
\Einheit=0.5cm
\[
\Pfad(-11,0),3433343334444343344344\endPfad
\SPfad(-11,0),1111111111111111111111\endSPfad
\DuennPunkt(-11,0)
\DuennPunkt(-10,1)
\DuennPunkt(-9,0)
\DuennPunkt(-8,1)
\DuennPunkt(-7,2)
\DuennPunkt(-6,3)
\DuennPunkt(-5,2)
\DuennPunkt(-4,3)
\DuennPunkt(-3,4)
\DuennPunkt(-2,5)
\DuennPunkt(-1,4)
\DuennPunkt(0,3)
\DuennPunkt(1,2)
\DuennPunkt(2,1)
\DuennPunkt(3,2)
\DuennPunkt(4,1)
\DuennPunkt(5,2)
\DuennPunkt(6,3)
\DuennPunkt(7,2)
\DuennPunkt(8,1)
\DuennPunkt(9,2)
\DuennPunkt(10,1)
\DuennPunkt(11,0)
\Label\u{ \textrm{ Dyck path}}(0,-1.0)
\Label\o{ \textrm{{\footnotesize \raisebox{2mm}{peak upstep} $\searrow$}} }(-9.1,2.6)
\]

\vspace*{5mm}

Next, an \emph{excedance} in a permutation $\pi$ on $[n]$ is a pair $(i,\pi(i))$ with $\pi(i)>i$; $i$ is the excedance \emph{location} and $\pi(i)$ is the excedance \emph{value}. Similarly, we have non-excedance locations and values.
Reifegerste \cite[p.\,761]{132diagram} observes that 321-avoiding permutations are characterized by the condition that the subwords formed by the excedance values and the non-excedance values are both increasing, and thus a 321-avoiding permutation is uniquely determined by its excedances, and, important for our purposes, also by its non-excedances. 

Now the Mansour-Deng-Du bijection from 321-avoiding permutations of $[n]$ to Dyck $n$-paths has a simple description as follows, with $n=11$ and 
\[
\left(\begin{array}{ccccccccccc}
1 & 2 & 3 & 4 & 5 & 6 & 7 & 8 & 9 & 10 & 11\\
1 & 4 & 5 & 2 & 6 & 9 & 3 & 7 & 11 & 8 & 10 
\end{array}\right)
\]
as a working example.

1. \ Extract the non-excedances: 
\[
\left(\begin{array}{cccccc}
1&4&7&8&10&11\\
1&2&3&7&8&10 
\end{array}\right)\, .
\]
 
2. \ Delete the last entry in the top row and the first entry in the bottom row (necessarily $n$ and 1 respectively), subtract 1 from each remaining entry in the bottom row, and align the rows:
\[
\left(\begin{array}{ccccc}
1&4&7&8&10\\
1&2&6&7&9
\end{array}\right)\, .
\]

This is the ascent-descent code of the desired Dyck path. Thus our working example corresponds to the Dyck path shown above.

Reifegerste's characterization above has an equivalent form: a permutation $\pi$ of $[n$] is 321-avoiding if and only if its excedance values are increasing and every \emph{non}-excedance value is a right to left minimum. (Note that an \emph{excedance} value $\pi(i)$ can never be a right to left minimum because there are too many entries after $\pi(i)$ for them all to be $>\pi(i)$.) This means that the non-excedance locations in a 321-avoiding permutation coincide with the locations of the peak upsteps (among all upsteps) in the corresponding Dyck path, and so one can verify that the mapping of Section 2 is indeed a bijection.

\end{document}